\documentclass[psfig,twocolumn,10pt]{asme2e}
\usepackage{amssymb}
\usepackage{amsfonts}
\usepackage{psfig}
\usepackage{graphicx}

\confshortname{IDETC/CIE 2013}
\conffullname{the ASME 2013 International Design Engineering Technical Conferences \&\\
              Computers and Information in Engineering Conference}
\confdate{12-15}
\confmonth{August}
\confyear{2013}
\confcity{Portland, OR}
\confcountry{USA}
\papernum{DETC2013-12151}
\input{tcilatex}
\begin{document}


\parindent 0pt \setcounter{topnumber}{9} \setcounter{bottomnumber}{9} %
\renewcommand{\textfraction}{0.00001}

\renewcommand {\floatpagefraction}{0.999} \renewcommand{\textfraction}{0.01} %
\renewcommand{\topfraction}{0.999} \renewcommand{\bottomfraction}{0.99} %
\renewcommand{\floatpagefraction}{0.99} \setcounter{totalnumber}{9}

\author{
	Andreas M\"uller$^1$, Zdravko Terze$^2$
	\affiliation{$^1$ Institute of Mechatronics, Reichenhainer Str. 88, 09126 Chemnitz, Germany, andreas.mueller@ifm-chemnitz.de\\
	$^2$ Faculty of Mechanical Eng. and Naval Arch., University of Zagreb, Croatia, Email: zdravko.terze@fsb.hr
	} 
\vspace{-3ex}
} 

\title{ 
\vspace{-4mm}
Is there an optimal choice of configuration space for Lie group integration schemes applied to constrained MBS?
\vspace{-2mm}
} 
\maketitle 

\begin{abstract}
{\it
Recently various numerical integration schemes have been proposed for
numerically simulating the dynamics of constrained multibody systems (MBS)
operating. These integration schemes operate directly on the MBS
configuration space considered as a Lie group. For discrete spatial
mechanical systems there are two Lie group that can be used as configuration
space: $SE\left( 3\right) $ and $SO\left( 3\right) \times \mathbb{R}^{3}$.
Since the performance of the numerical integration scheme clearly depends on
the underlying configuration space it is important to analyze the effect of
using either variant. For constrained MBS a crucial aspect is the constraint
satisfaction. In this paper the constraint violation observed for the two
variants are investigated. It is concluded that the $SE\left( 3\right) $
formulation outperforms the $SO\left( 3\right) \times \mathbb{R}^{3}$
formulation if the absolute motions of the rigid bodies, as part of a
constrained MBS, belong to a motion subgroup. In all other cases both
formulations are equivalent. In the latter cases the $SO\left( 3\right)
\times \mathbb{R}^{3}$ formulation should be used since the $SE\left(
3\right) $ formulation is numerically more complex, however.
}

{\it Keywords--} {Constrained multibody systems, Lie group integration, screw systems}

\end{abstract}

\section{Introduction}

Lie group integration schemes for MBS commonly rest on $SO\left( 3\right)
\times \mathbb{R}^{3}$ as configuration space manifold. This configuration
space cannot capture the intrinsic geometry of rigid body motions since it
does not represent proper screw motions. Moreover the general motion of a
rigid body is a screw motion. This applies to unconstrained as well as
constrained rigid bodies, and the reconstruction of finite motions within
numerical time stepping schemes shall take this into account. Along this
line the Lie group $SE\left( 3\right) $ of proper rigid body motions was
recently used as configuration space \cite{BruelsCardona2010},\cite%
{BruelsCardonaArnold2012},\cite{Terze2011}. It turned out that the proper
rigid body motion group does outperform the standard $SO\left( 3\right)
\times \mathbb{R}^{3}$ formulation for a rigid bodies constrained to move
relative to a stationary reference (e.g. heavy top), The obvious question is
whether this statement applies to general constrained MBS, and whether there
is an optimal choice for a given MBS that leads to the best numerical
performance. In this paper two Lie group formulations for constrained MBS
are compared and the effect of using either configuration space is analyzed
for several examples.

In Lie group setting the dynamics of a constrained multibody system (MBS) is
governed by the constrained Boltzmann-Hamel equations%
\begin{eqnarray}
\mathbf{M}\left( q\right) \dot{\mathbf{V}}+\mathbf{J}^{T}\mathbf{\lambda }
&=&\mathbf{Q}\left( q,\mathbf{V},t\right) {\tag{1a}}  \label{dyn1} \\
\dot{q} &=&q\mathbf{V}{\tag{1b}} \\
\mathbf{g}\left( q\right) &=&\mathbf{0}{\tag{1c}}\setcounter{equation}{2}
\end{eqnarray}%
where $q\in G$ represents the MBS configuration and $G$ is the configuration
space Lie group. This is an index 3 DAE system on the Lie group $G$. The
system (\ref{dyn1}a) represents the motion equations of the MBS subjected to
the geometric constraints (\ref{dyn1}c) that are complemented by the \textit{%
kinematic reconstruction equations} (\ref{dyn1}b). That is, integration of (%
\ref{dyn1}b) yields the motion $\mathbf{q}\left( t\right) $ of the MBS
corresponding to the MBS velocity $\mathbf{V}\in \mathfrak{g}$, with $%
\mathfrak{g}$ being the Lie algebra of $G$. In order to apply Lie group ODE
integration schemes the first step is the reformulation of (\ref{dyn1}) as
ODE on the state space $S:=G\times \mathfrak{g}$. This is achieved with the
widely used index 1 formulation%
\begin{equation}
\left( 
\begin{array}{cc}
\mathbf{M} & \mathbf{J}^{T} \\ 
\mathbf{J} & \mathbf{0}%
\end{array}%
\right) \left( 
\begin{array}{c}
\dot{\mathbf{V}} \\ 
\mathbf{\lambda }%
\end{array}%
\right) =\left( 
\begin{array}{c}
\mathbf{Q} \\ 
\mathbf{\eta }%
\end{array}%
\right)  \label{dyn4}
\end{equation}%
using the acceleration constraints $\mathbf{J}\left( q\right) \cdot \dot{%
\mathbf{V}}=\mathbf{\eta }\left( q,\mathbf{V}\right) $. This system replaces
the dynamic equations (\ref{dyn1}a) when subject to the holonomic
constraints (\ref{dyn1}c) since, for a given state $X=\left( q,\mathbf{V}%
\right) \in S$, the system (\ref{dyn4}) and thus (\ref{dyn1}) can be solved
for $\dot{\mathbf{V}}$ consistent with the acceleration constraints. The
index reduction is achieved on the expense of numerical drifts of the
constraints, however \cite{Blajer2011}. The overal ODE system is obtained by
complimenting (\ref{dyn4}) with the kinematic equations (\ref{dyn1}b).

As far as the MBS model is concerned the numerical performance, and
eventually the accuracy of any integration scheme whatsoever, are determined
by 1) the choice of generalized coordinates and 2) by how generalized
velocities are introduced. The first issue has to do with finding a proper
chart on the configuration space, whereas the latter concerns the relation
of velocities and time derivatives of the configurations, i.e. the relation (%
\ref{dyn1}b). The best solution for the first issue is to avoid local
coordinates at all. The Lie group concept provides such a geometric vehicle
that allows for coordinate-free modeling of frame transformations, where
either $SO\left( 3\right) \times \mathbb{R}^{3}$ or $SE\left( 3\right) $ can
be used. The second issue concerns the particular configuration space Lie
group, noting that $q\in G$ and $\mathbf{V}\in \mathfrak{g}$.

The integration method considered in this paper is the Munthe-Kaas method.
Since the kinematic reconstruction is inherent to the model the
considerations shall apply to Lie group integration schemes in general.

\section{Two State Spaces for Rigid Bodies}

The configuration of a rigid body, with respect to a space-fixed inertial
reference frame, is described by the position vector $\mathbf{r}\in \mathbb{R%
}^{3}$ of the origin of a body-fixed reference frame and its rotation matrix 
$\mathbf{R}\in SO\left( 3\right) $, summarized by the pair $C=\left( \mathbf{%
R},\mathbf{r}\right) $. A rigid body motion is thus a curve $C\left(
t\right) $. The crucial point is to assign the Lie group $C$ is living in.
The state space is then the product of this Lie group and its Lie algebra.

\subsection{Group of Proper Rigid Body Motions $SE\left( 3\right) $}

$SE\left( 3\right) $ represents frame transformations, i.e. the combination
of two successive rigid-body comfigurations are given by $C_{2}\cdot
C_{1}=\left( \mathbf{R}_{2}\mathbf{R}_{1},\mathbf{r}_{2}+\mathbf{R}_{2}%
\mathbf{r}_{1}\right) $. This multiplication rule indicates that $SE\left(
3\right) :=SO\left( 3\right) \ltimes \mathbb{R}^{3}$ is the 6-dimensional
semidirect product group of the rotation group $SO\left( 3\right) $ and the
translation group, represented as $\mathbb{R}^{3}$. Rigid body
configurations can be represented in matrix form as%
\begin{equation}
\mathbf{C}=\left( 
\begin{array}{cc}
\mathbf{R} & \mathbf{r} \\ 
\mathbf{0} & 1%
\end{array}%
\right)  \label{C}
\end{equation}%
which admits representing the group multiplication as matrix multiplication%
\begin{equation}
\mathbf{C}_{2}\mathbf{C}_{1}=\left( 
\begin{array}{cc}
\mathbf{R}_{2}\mathbf{R}_{1} & \mathbf{r}_{2}+\mathbf{R}_{2}\mathbf{r}_{1}
\\ 
\mathbf{0} & 1%
\end{array}%
\right) .  \label{MultSE3}
\end{equation}

A generic motion of a rigid body is a screw motion, i.e. an interconnected
rotation and translation along a screw axis \cite{BottemaRoth1979}.

The velocity corresponding to the screw motion of a rigid body is a twist $%
\mathbf{V}=\left( \mathbf{\omega },\mathbf{v}\right) \in \mathbb{R}^{6}$
with angular velocity $\mathbf{\omega }$ and linear velocity vector $\mathbf{%
v}$. The \textit{body-fixed twist} of a rigid body motion $\mathbf{C}\left(
t\right) $ is determined as 
\begin{equation}
\widehat{\mathbf{V}}:=\mathbf{C}^{-1}\dot{\mathbf{C}}\text{ \ \ with \ \ }%
\widehat{\mathbf{V}}=\left( 
\begin{array}{cc}
\widehat{\mathbf{\omega }} & \mathbf{v} \\ 
\mathbf{0} & 0%
\end{array}%
\right) \in se\left( 3\right)  \label{Vhat}
\end{equation}%
where $se\left( 3\right) $ is the Lie algebra of $SE(3)$. The assignment (%
\ref{Vhat}) is a one-one correpsondance of twist coordinate vectors and $%
se\left( 3\right) $-matrices. $\widehat{\mathbf{\omega }}:=\mathbf{R}^{T}%
\dot{\mathbf{R}}\in so\left( 3\right) $ is the skew symmetric cross product
matrix associated to the vector $\mathbf{\omega }$. Via this isomorphism the
Lie bracket on $se\left( 3\right) $, in vector representation, is given the
screw product $\left[ \mathbf{V}_{1},\mathbf{V}_{2}\right] =\left( \mathbf{%
\omega }_{1}\times \mathbf{\omega }_{2},\mathbf{\omega }_{1}\times \mathbf{v}%
_{2}-\mathbf{\omega }_{2}\times \mathbf{v}_{1}\right) $ \cite%
{BottemaRoth1979}. For convenience $\mathbf{\omega }\in so\left( 3\right) $
is written for a vector $\mathbf{\omega }\in \mathbb{R}^{3}$ to indicate
that isomorphism of $so\left( 3\right) $ and $\mathbb{R}^{3}$ equipped with
the cross product as Lie bracket.

For any screw $\mathbf{X}=\left( \mathbf{\omega },\mathbf{v}\right) \in 
\mathbb{R}^{6}$ with screw axis parallel to $\mathbf{\omega }$ the linear
part can be expressed with a position vector $\mathbf{r}$ on the screw axis
as $\mathbf{v}=\mathbf{r}\times \mathbf{\omega }+h\mathbf{\omega }$. The
screw $\mathbf{X}$ describes an instantaneous screw motion, i.e. a rotation
about the axis $\mathbf{\omega }$ together with a translation along this
axis, where $h:=\mathbf{\omega }\cdot \mathbf{v}/\left\Vert \mathbf{\omega }%
\right\Vert ^{2}$ is the pitch of the screw. If $h=0$, then $\mathbf{X}$ are
simply the Pl\"{u}cker coordinates of a line parallel to the screw axis.

The Lie bracket can be expressed by a linear operation on screw coordinate
vectors as $[{\mathbf{V}}_{1},{\mathbf{V}}_{2}]=\mathbf{ad}_{{\mathbf{V}}%
_{1}}{\mathbf{V}}_{2}$ given by the matrix 
\begin{equation}
\mathbf{ad}_{\mathbf{V}}=\left( 
\begin{array}{cc}
\widehat{\mathbf{\omega }} & 0 \\ 
\widehat{\mathbf{v}} & \widehat{\mathbf{\omega }}%
\end{array}%
\right) .  \label{adse3}
\end{equation}

The connection between an infinitesimal screw motion $\mathbf{X}\left(
t\right) $ and the correponding finite screw motion is given by the
exponential mapping on $SE\left( 3\right) $, which reads explicitly%
\begin{equation}
\mathbf{X}=\left( \mathbf{\omega },\mathbf{v}\right) \longmapsto \exp 
\widehat{\mathbf{X}}=\left( 
\begin{array}{cc}
\exp \widehat{\mathbf{\omega }} & \;\;\;\frac{1}{\left\Vert \mathbf{\omega }%
\right\Vert ^{2}}\left( I-\exp \widehat{\mathbf{\omega }}\right) \left( 
\mathbf{\omega }\times \mathbf{v}\right) +h\,\mathbf{\omega }%
\vspace{2mm}
\\ 
0 & 1%
\end{array}%
\right)  \label{expX}
\end{equation}%
where%
\vspace{-2ex}
\begin{equation}
\exp \widehat{\mathbf{\omega }}=\mathbf{I}+\frac{\sin \left\Vert \mathbf{%
\omega }\right\Vert }{\left\Vert \mathbf{\omega }\right\Vert }\widehat{%
\mathbf{\omega }}+\frac{1-\cos \left\Vert \mathbf{\omega }\right\Vert }{%
\left\Vert \mathbf{\omega }\right\Vert }\widehat{\mathbf{\omega }}^{2}
\label{expSO3}
\end{equation}%
is the exponential mapping on $SO\left( 3\right) $. Important for the Lie
group integration scheme is the differential of the exponential mapping, $%
\mathrm{dexp}:se\left( 3\right) \times se\left( 3\right) \rightarrow
se\left( 3\right) $ that can be introduced as $\mathrm{dexp}_{\widehat{%
\mathbf{X}}}\dot{\widehat{\mathbf{X}}}=\dot{\mathbf{C}}\mathbf{C}^{-1}$,
with $\mathbf{C}=\exp \widehat{\mathbf{X}}$. Its signinifcance becomes clear
by replacing $\mathbf{X}$ with $-\mathbf{X}$, which allows expressing the
body-fixed twist as $\widehat{\mathbf{V}}=\mathrm{dexp}_{-\widehat{\mathbf{X}%
}}\dot{\widehat{\mathbf{X}}}$. The numerical ODE integration methods require
the inverse of the dexp map. In vector representation of twists the inverse
of dexp on $SE\left( 3\right) $ is \cite{SeligBook}%
\begin{eqnarray}
\mathrm{dexp}_{\widehat{\mathbf{X}}}^{-1} &=&\mathbf{I}-\frac{1}{2}\mathbf{ad%
}_{\mathbf{X}}+\left( \frac{2}{\left\Vert \mathbf{\omega }\right\Vert ^{2}}+%
\frac{\left\Vert \mathbf{\omega }\right\Vert +3\sin \left\Vert \mathbf{%
\omega }\right\Vert }{4\left\Vert \mathbf{\omega }\right\Vert \left( \cos
\left\Vert \mathbf{\omega }\right\Vert -1\right) }\right) \mathbf{ad}_{%
\mathbf{X}}^{2}  \nonumber \\
&&\hspace{1.5ex}+\left( \frac{1}{\left\Vert \mathbf{\omega }\right\Vert ^{4}}%
+\frac{\left\Vert \mathbf{\omega }\right\Vert +\sin \left\Vert \mathbf{%
\omega }\right\Vert }{4\left\Vert \mathbf{\omega }\right\Vert ^{3}\left(
\cos \left\Vert \mathbf{\omega }\right\Vert -1\right) }\right) \mathbf{ad}_{%
\mathbf{X}}^{4}  \label{dexpInvSE3}
\end{eqnarray}%
with $\mathbf{X}=\left( \mathbf{\omega },\mathbf{v}\right) $. In vector
representation of $so(3)$ the inverse of the differential of the exp mapping
(\ref{expSO3}) for $\mathbf{R}=\exp \mathbf{\xi }$ is given as matrix \cite%
{BulloMurray1995} 
\begin{equation}
\mathrm{dexp}_{\mathbf{\xi }}^{-1}=\mathbf{I}-\frac{1}{2}\widehat{\mathbf{%
\xi }}+\left( 1-\frac{\left\Vert \mathbf{\xi }\right\Vert }{2}\cot \frac{%
\left\Vert \mathbf{\xi }\right\Vert }{2}\right) \frac{\widehat{\mathbf{\xi }}%
^{2}}{\left\Vert \mathbf{\xi }\right\Vert ^{2}}.  \label{dexpInvSO3}
\end{equation}

The state of a single unconstrained rigid body is represented by the couple $%
\left( \mathbf{C},\mathbf{V}\right) \in SE\left( 3\right) \times se\left(
3\right) $. The acceleration of the body is the time derivative $\dot{%
\mathbf{V}}\in \mathbb{R}^{6}$. Making use of (\ref{Vhat}) the time
derivative of $\mathbf{C}$ can be identified with $\mathbf{V}$. The time
derivative of the rigid body state $\left( \mathbf{C},\mathbf{V}\right) $ is
thus isomorphic to $(\widehat{\mathbf{V}},\dot{\mathbf{V}})\in se\left(
3\right) \times \mathbb{R}^{6}$.

Consequently, in $SE\left( 3\right) $ representation, the state space of a
rigid body is the Lie group $SE\left( 3\right) \times \mathbb{R}^{6}$ with
Lie algebra $se\left( 3\right) \times \mathbb{R}^{6}$. The multiplication on
this rigid body state space is%
\begin{equation}
\left( \mathbf{C}_{1},\mathbf{V}_{1}\right) \cdot \left( \mathbf{C}_{2},%
\mathbf{V}_{2}\right) =\left( \mathbf{C}_{1}\mathbf{C}_{2},\mathbf{V}_{1}+%
\mathbf{V}_{2}\right) .
\end{equation}%
Being a Lie group the state space possesses an exponential mapping given by%
\begin{equation}
\begin{array}{cccc}
\exp :\ \  & se\left( 3\right) \times \mathbb{R}^{6} & \rightarrow & 
SE\left( 3\right) \times \mathbb{R}^{6} \\ 
& (\widehat{\mathbf{V}},\mathbf{A)} & \mapsto & (\exp \widehat{\mathbf{V}},%
\mathbf{A)}%
\end{array}
\label{expSBody}
\end{equation}%
with the exponential (\ref{expX}). The rigid body state can thus be
reconstructed from its time derivative via this exponential mapping. With
the Lie bracket on the algebra $se\left( 3\right) \times \mathbb{R}^{6}$%
\begin{equation}
\lbrack (\widehat{\mathbf{V}}_{1},\mathbf{A}_{1}),(\widehat{\mathbf{V}}_{2},%
\mathbf{A}_{2})]=([\widehat{\mathbf{V}}_{1},\widehat{\mathbf{V}}_{2}],%
\mathbf{0})
\end{equation}%
the differential of the exponential mapping is%
\begin{equation}
\mathrm{dexp}_{(\widehat{\mathbf{V}}_{1},\mathbf{A}_{1}\mathbf{)}}(\widehat{%
\mathbf{V}}_{2},\mathbf{A}_{2}\mathbf{)}=(\mathrm{dexp}_{\widehat{\mathbf{V}}%
_{1}}\widehat{\mathbf{V}}_{2},\mathbf{A}_{2}).
\end{equation}

\subsection{Direct Product Group $SO\left( 3\right) \times \mathbb{R}^{3}$}

Neglecting the interrelation of rotations and translations the
multiplication is%
\begin{equation}
C_{1}\cdot C_{2}=(\mathbf{R}_{1}\mathbf{R}_{2},\mathbf{r}_{1}+\mathbf{r}_{2})
\label{MultDirectProduct}
\end{equation}%
which indicates that $C=\left( \mathbf{R},\mathbf{r}\right) \in SO\left(
3\right) \times \mathbb{R}^{3}$. The direct product $SO\left( 3\right)
\times \mathbb{R}^{3}$ is commonly used as rigid body configuration space
for Lie group methods \cite%
{BruelsCardona2010,BruelsCardonaArnold2012,CelledoniOwren1999,Krysl,Terze2011}%
. Apparently this multiplication does not represent frame transformations.
The inverse element is $C^{-1}=(\mathbf{R}^{T},-\mathbf{r})$.

The Lie algebra of the direct product $SO\left( 3\right) \times \mathbb{R}%
^{3}$ is $so\left( 3\right) \times \mathbb{R}^{3}$ with Lie bracket%
\begin{equation}
\left[ \mathbf{X}_{1},\mathbf{X}_{2}\right] =(\mathbf{\omega }_{1}\times 
\mathbf{\omega }_{2},\mathbf{0}).  \label{LieBracketdirect}
\end{equation}%
The exponential mapping on the direct product group is%
\begin{equation}
\mathbf{X}=\left( \mathbf{\omega },\mathbf{v}\right) \longmapsto \exp 
\mathbf{X}=\left( \exp \widehat{\mathbf{\omega }},\mathbf{v}\right)
\label{expSO3xR3}
\end{equation}%
with the exponential mapping (\ref{expSO3}) on $SO\left( 3\right) $. The
dexp mapping is accordingly%
\begin{equation}
\mathrm{dexp}_{\left( \mathbf{\xi ,u}\right) }\left( \mathbf{\eta },\mathbf{v%
}\right) =(\mathrm{dexp}_{\mathbf{\xi }}\mathbf{\eta },\mathbf{v}),
\label{dexpSO3xR3}
\end{equation}%
with dexp mapping on $SO\left( 3\right) $. Its inverse is readily $\mathrm{%
dexp}_{\left( \mathbf{\xi ,u}\right) }^{-1}\left( \mathbf{\eta },\mathbf{v}%
\right) =(\mathrm{dexp}_{\mathbf{\xi }}^{-1}\mathbf{\eta },\mathbf{v})$.

The velocity of a rigid body is, with configuration $C\in SO\left( 3\right)
\times \mathbb{R}^{3}$, given as%
\begin{equation}
\left( \widehat{\mathbf{\omega }},\mathbf{v}^{s}\right) =(\mathbf{R}^{T}\dot{%
\mathbf{R}},\dot{\mathbf{r}}):=C^{-1}\dot{C}\in so\left( 3\right) \times 
\mathbb{R}^{3}  \label{hybridvel}
\end{equation}%
and in vector notation denoted $\mathbf{V}=\left( \mathbf{\omega },\mathbf{v}%
^{s}\right) \in \mathbb{R}^{3}\times \mathbb{R}^{3}$, with $\mathbf{v}^{s}=%
\dot{\mathbf{r}}$. This velocity couple is clearly not a proper twist but
contains a mix of body-fixed angular velocity $\mathbf{\omega }$ and spatial
linear velocity $\mathbf{v}^{s}$. It is therefore called \textit{hybrid
representation} of rigid body velocities \cite{Bruyninckx1996},\cite{murray}%
. It is frequently used for expressing the Newton-Euler equations. Even
though, angular and linear velocities are treated independently, and this
definition does not reflect the intrinsic characteristics of screw motions.

In hybrid representation the state of a rigid body is represented by $\left(
C,\mathbf{V}\right) \in SO\left( 3\right) \times \mathbb{R}^{3}\times
so\left( 3\right) \times \mathbb{R}^{3}$. This is a Lie group with algebra $%
so\left( 3\right) \times \mathbb{R}^{3}\times \mathbb{R}^{6}$.
Multiplication is defined as%
\begin{equation}
\left( C_{1},\mathbf{V}_{1}\right) \cdot \left( C_{2},\mathbf{V}_{2}\right)
=\left( C_{1}\cdot C_{2},\mathbf{V}_{1}+\mathbf{V}_{2}\right) ,
\end{equation}%
and the exponential mapping is, with (\ref{expSO3xR3}),%
\begin{equation}
\begin{array}{cccc}
\exp :\ \  & so\left( 3\right) \times \mathbb{R}^{3}\times \mathbb{R}^{6} & 
\rightarrow & SO\left( 3\right) \times \mathbb{R}^{3}\times \mathbb{R}%
^{3}\times \mathbb{R}^{3} \\ 
& (\mathbf{V},\mathbf{A)} & \mapsto & (\exp \mathbf{V},\mathbf{A)\ .}%
\end{array}%
\end{equation}%
The Lie bracket is on this algebra is%
\begin{equation}
\lbrack (\mathbf{V}_{1},\mathbf{A}_{1}),(\mathbf{V}_{2},\mathbf{A}_{2})]=([%
\mathbf{V}_{1},\mathbf{V}_{2}],\mathbf{0}).
\end{equation}%
The differential of the exponential mapping is, with (\ref{dexpSO3xR3}),
given by%
\begin{equation}
\mathrm{dexp}_{(\mathbf{V}_{1},\mathbf{A}_{1})}(\mathbf{V}_{2},\mathbf{A}_{2}%
\mathbf{)}=(\mathrm{dexp}_{\mathbf{V}_{1}}\mathbf{V}_{2},\mathbf{A}_{2}).
\end{equation}

\section{State Space of MBS}

\subsection{Group of Proper Rigid Body Motions $SE\left( 3\right) $}

The configuration of an MBS consisting of $n$ rigid bodies is represented by 
$q=(\mathbf{C}_{1},\ldots ,\mathbf{C}_{n})\in G^{\ltimes }$, where%
\begin{equation}
G^{\ltimes }:=SE\left( 3\right) ^{n}
\end{equation}%
is the $6n$-dimensional configuration space Lie group. This is a
coordinate-free representation of MBS configurations. The multiplication on $%
G^{\ltimes }$ is understood componentwise, and thus $q^{-1}=(\mathbf{C}%
_{1}^{-1},\ldots ,\mathbf{C}_{n}^{-1})$. The MBS velocities are represented
by $\mathbf{V}=(\mathbf{V}_{1},\ldots ,\mathbf{V}_{n})\in \left( \mathbb{R}%
^{6}\right) ^{n}$. The body-fixed velocities are determined by $\widehat{%
\mathbf{V}}=q^{-1}\dot{q}$ denoting $\widehat{\mathbf{V}}=(\widehat{\mathbf{V%
}}_{1},\ldots ,\widehat{\mathbf{V}}_{n})$. The MBS state space is thus the $%
12\cdot n$-dimensional Lie group%
\begin{equation}
S^{\ltimes }:=SE\left( 3\right) ^{n}\times \left( \mathbb{R}^{6}\right) ^{n}
\label{StateSpace}
\end{equation}%
and the MBS state is $X=\left( q,\mathbf{V}\right) =(\mathbf{C}_{1},\ldots ,%
\mathbf{C}_{n},\mathbf{V}_{1},\ldots ,\mathbf{V}_{n})\in S^{\ltimes }$. The
multiplication is $X^{\prime }\cdot X^{\prime \prime }=\left( \mathbf{C}%
_{1}^{\prime }\mathbf{C}_{1}^{\prime \prime },\ldots ,\mathbf{C}_{n}^{\prime
}\mathbf{C}_{n}^{\prime \prime },\mathbf{V}_{1}^{\prime }+\mathbf{V}%
_{1}^{\prime \prime },\ldots ,\mathbf{V}_{n}^{\prime }+\mathbf{V}%
_{n}^{\prime \prime }\right) $.

The corresponding Lie algebra is%
\begin{equation}
\mathfrak{s}^{\ltimes }:=se\left( 3\right) ^{n}\times \left( \mathbb{R}%
^{6}\right) ^{n},
\end{equation}%
with $x=(\mathbf{V}_{1},\ldots ,\mathbf{V}_{n},\mathbf{A}_{1},\ldots ,%
\mathbf{A}_{n})\in \mathfrak{s}^{\ltimes }$. The exponential mapping on the
state space is%
\begin{equation}
\exp x=\left( \exp \mathbf{V}_{1},\ldots ,\exp \mathbf{V}_{n},\mathbf{A}%
_{1},\ldots ,\mathbf{A}_{n}\right) \in S^{\ltimes }  \label{expS}
\end{equation}%
with (\ref{expSBody}) with differential $\mathrm{dexp}:\mathfrak{s}^{\ltimes
}\times \mathfrak{s}^{\ltimes }\rightarrow \mathfrak{s}^{\ltimes }$ 
\begin{equation}
\mathrm{dexp}_{x^{\prime }}x^{\prime \prime }=(\mathrm{dexp}_{\mathbf{V}%
_{1}^{\prime }}\mathbf{V}_{1}^{\prime \prime },\ldots ,\mathrm{dexp}_{%
\mathbf{V}_{n}^{\prime }}\mathbf{V}_{n}^{\prime \prime },\mathbf{A}%
_{1}^{\prime \prime },\ldots ,\mathbf{A}_{n}^{\prime \prime }).
\label{dexpS}
\end{equation}

\subsection{Direct Product Group $SO\left( 3\right) \times \mathbb{R}^{3}$}

When the direct product group is used the MBS configuration $q=(C_{1},\ldots
,C_{n})\in G^{\times }$ belongs to the $6n$-dimensional Lie group%
\begin{equation}
G^{\times }:=\left( SO\left( 3\right) \times \mathbb{R}^{3}\right) ^{n}
\end{equation}%
and possess the inverse $q^{-1}=(C_{1}^{-1},\ldots ,C_{n}^{-1})$. The hybrid
velocity of the MBS is $q^{-1}\dot{q}=\left( \left( \widehat{\mathbf{\omega }%
}_{1},\mathbf{v}_{1}^{s}\right) ,\ldots ,\left( \widehat{\mathbf{\omega }}%
_{n},\mathbf{v}_{n}^{s}\right) \right) $, and written as vector $V=\left( 
\mathbf{V}_{1},\ldots ,\mathbf{V}_{n}\right) \in \left( \mathbb{R}%
^{6}\right) ^{n}$ with $\mathbf{V}_{i}=\left( \mathbf{\omega }_{i},\mathbf{v}%
_{i}^{s}\right) $. Therewith the MBS state space is%
\begin{equation}
S^{\times }:=\left( SO\left( 3\right) \times \mathbb{R}^{3}\right)
^{n}\times \left( \mathbb{R}^{6}\right) ^{n}
\end{equation}%
with state vector $X=\left( q,\mathbf{V}\right) =(C_{1},\ldots ,C_{n}\mathbf{%
,V}_{1},\ldots ,\mathbf{V}_{n})\in S^{\times }$. This is a $12\cdot n$%
-dimensional Lie group with multiplication $X^{\prime }\cdot X^{\prime
\prime }=(C_{1}^{\prime }\cdot C_{1}^{^{\prime \prime }},\ldots
,C_{n}^{\prime }\cdot C_{n}^{^{\prime \prime }}\mathbf{,V}_{1}^{\prime }+%
\mathbf{V}_{1}^{\prime \prime },\ldots ,\mathbf{V}_{n}^{\prime }+\mathbf{V}%
_{n}^{\prime \prime })$. The Lie algebra of $S^{\times }$ is%
\begin{equation}
\mathfrak{s}^{\times }:=\left( so\left( 3\right) \times \mathbb{R}%
^{3}\right) ^{n}\times \left( \mathbb{R}^{6}\right) ^{n},
\end{equation}%
with elements $x=(\mathbf{V}_{1},\ldots ,\mathbf{V}_{n},\mathbf{A}%
_{1},\ldots ,\mathbf{A}_{n})\in \mathfrak{s}^{\times }$ where $\mathbf{A}%
_{i}=\left( \mathbf{\alpha }_{i},\mathbf{a}_{i}^{s}\right) $ represents the
body-fixed angular and spatial linear acceleration. The exponential mapping
on the ambient state space is%
\begin{equation}
\exp x=\left( \exp \mathbf{V}_{1},\ldots ,\exp \mathbf{V}_{n},\mathbf{A}%
_{1},\ldots ,\mathbf{A}_{n}\right) \in S^{\times }
\end{equation}%
with (\ref{expSO3xR3}). The differential $\mathrm{dexp}:\mathfrak{s}^{\times
}\times \mathfrak{s}^{\times }\rightarrow \mathfrak{s}^{\times }$ is, with (%
\ref{dexpSO3xR3}),%
\begin{equation}
\mathrm{dexp}_{x^{\prime }}x^{\prime \prime }=(\mathrm{dexp}_{\mathbf{V}%
_{1}^{\prime }}\mathbf{V}_{1}^{\prime \prime },\ldots ,\mathrm{dexp}_{%
\mathbf{V}_{n}^{\prime }}\mathbf{V}_{n}^{\prime \prime },\mathbf{A}%
_{1},\ldots ,\mathbf{A}_{n}).
\end{equation}

\section{Munthe-Kaas Method for Constrained MBS}

The Munthe-Kaas (MK) method \cite{muntekaas1998,muntekaas1999,Hairer2006} is
a widely used integration scheme for ODE on Lie group that was applied to
rigid body dynamics such as \cite{CelledoniOwren1999}. Its appeal stems from
its construction since it is a direct extension of the Runge-Kutta method.
In order to apply MK a scheme the system equations must be expressed in the
form of a first-order ODE on the state space%
\begin{equation}
\dot{X}=XF\left( t,X\right)  \label{ODEleft}
\end{equation}%
with a mapping $F:\mathbb{R}\times S\rightarrow \mathfrak{s}$. This is
achieved with help of the index 1 system (\ref{dyn4}). In order to solve for 
$\dot{X}$, at a given state $X=\left( q,\mathbf{V}\right) \in S$, the system
(\ref{dyn4}) must be solved for $\dot{\mathbf{V}}$. It remains to evaluate (%
\ref{dyn1}b) for $\dot{q}$. By introducing the mapping $F\left( t,X\right) =(%
\mathbf{V},\dot{\mathbf{V}})$, which includes solving for $\dot{\mathbf{V}}$%
, the system (\ref{dyn1}) is equivalent to (\ref{ODEleft}) since $XF\left(
t,X\right) =(q\mathbf{V},\dot{\mathbf{V}})$. The equations (\ref{ODEleft})
can be regarded as the Boltzmann-Hamel equations on the state space Lie
group when its algebra is defined by left trivialization. Evaluation of $%
XF\left( t,X\right) $ thus amounts to solving (\ref{dyn4}) for $\dot{\mathbf{%
V}}$ and evaluating (\ref{dyn1}b). Since both, the rigid body twists (\ref%
{Vhat}) and the hybrid velocities (\ref{hybridvel}) are defined by left
translation this applies to the $SE\left( 3\right) $ and the $SO\left(
3\right) \times \mathbb{R}^{3}$ formulation as e.g. in \cite%
{BruelsCardona2010,BruelsCardonaArnold2012,Terze2011}.

In the MK method solutions are sought of the form $X\left( t\right)
=X_{0}\exp \Phi \left( t\right) $. This allows replacing the original system
(\ref{ODEleft}) at the integration step $i$ by the system%
\begin{equation}
\dot{\Phi}^{\left( i\right) }=\mathrm{dexp}_{-\Phi ^{\left( i\right)
}}^{-1}F(t,X_{i-1}\exp \Phi ^{\left( i\right) }),\ t\in \lbrack
t_{i-1},t_{i}],\text{ with }\Phi ^{\left( i\right) }\left( t_{i-1}\right) =0
\label{subst1}
\end{equation}%
with initial condition $X_{i-1}$. Notice the negative sign of $\Phi ^{\left(
i\right) }$ in dexp, which is different from the original MK version.
Originally the MK method is derived for right invariant systems, i.e. ODE
systems of the form $\dot{X}=F\left( t,X\right) X$. Numerically solving (\ref%
{subst1}) yields a solution $\Phi ^{\left( i\right) }\left( t_{i}\right) $,
and thus a solution $X_{i}:=X_{i-1}\exp \Phi ^{\left( i\right) }\left(
t_{i}\right) $ of (\ref{ODEleft}). The $\Phi ^{\left( i\right) }$ represent
local coordinates on the state space defined in a neighborhood of $X_{i-1}$.
The system (\ref{subst1}) is solved by an $s$-stage RK method. This gives
rise the corresponding $s$-stage MK scheme at time step $i$ 
\begin{eqnarray}
X_{i}{:=} &&X_{i-1}\exp \Phi ^{\left( i\right) },\ \Phi ^{\left( i\right) }{%
:=}h\sum_{j=1}^{s}b_{j}k_{j}  \label{MK1} \\
k_{j}{:=} &&\mathrm{dexp}_{-\Psi _{j}}^{-1}F\left(
t_{i-1}+c_{j}h,X_{i-1}\exp \Psi _{j}\right)  \nonumber \\
\Psi _{j}\,{:=} &&h\sum_{l=1}^{j-1}a_{jl}k_{l},\ \Psi _{1}=0,  \nonumber
\end{eqnarray}%
where $a_{jl},b_{j}$, and $c_{j}$ are the Butcher coefficients of the $s$%
-stage RK method, and $k_{j},\Psi _{j}\in \mathfrak{s}$.

It is well-known that the index 1 formulation suffers from constraint
violations due to numerical drifts introduced by the numerical update
scheme. This is indeed carries over to the introduced Lie group formulation,
and the established constraint stabilization methods can be applied. The
paper \cite{Blajer2011} provides a good overview of constraint stabilization
methods.

\section{Examples}

\subsection{Heavy Top in Gravity Field}

As first example a heavy top is considered, i.e. a single rigid body
constrained to rotate about a space-fixed pivot point.

\begin{figure}[h]
\centerline{
\includegraphics[width=4.5cm]{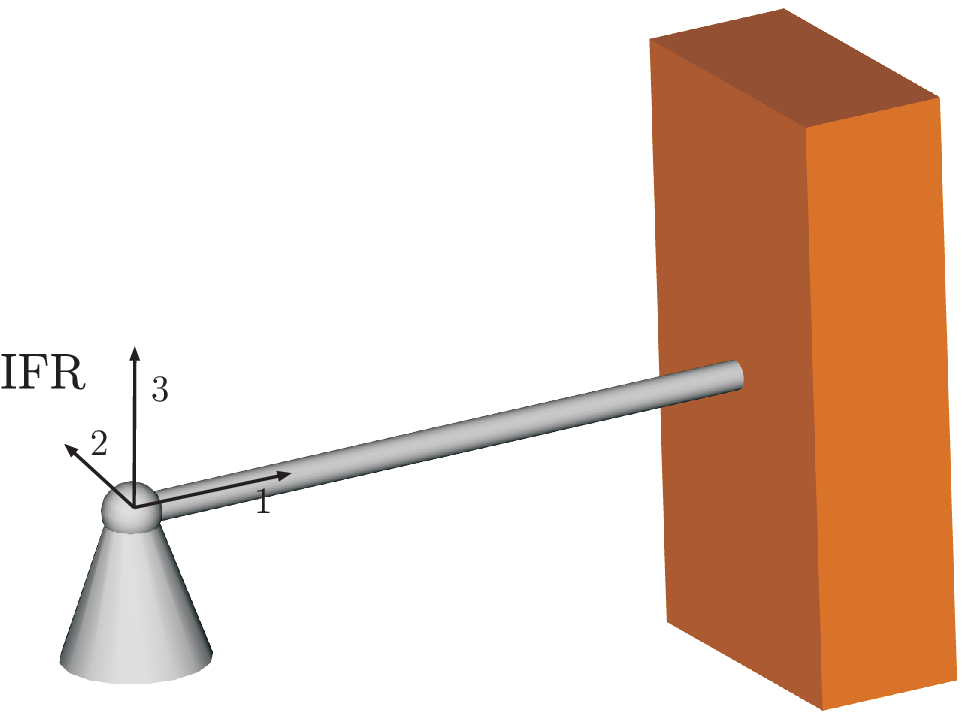}
} \vspace{1.5ex} \vspace{-3ex}
\caption{Model of a heavy top.}
\label{figHaevyTop}
\end{figure}

The model concists of a rectangular solid box with side lengths $0.1\times
0.2\times 0.4$\thinspace m connected to the ground by a speherical joint as
shown in figure \ref{figHaevyTop}a). A body-fixed reference frame (RFR) is
attached at the COM. In the shown reference configuration the RFR is aligned
to the space-fixed inertia frame (IFR). Assuming aluminium material the body
mass is $m=21.6$~kg, and its inertia tensor w.r.t. the RFR is ${\Theta }_{0}=%
\mathrm{diag~}\left( 0.36,{0.306,0.09}\right) ~$kg~m$^{2}$. The position
vector of the COM measured from the pivot point expressed in the body-fixed
reference frame is denoted with $\mathbf{r}_{0}=(0.5,0,0)$\thinspace m. The
configuration of the reference frame is represented by $C=\left( \mathbf{R},%
\mathbf{r}^{s}\right) $. \vspace{-3ex}

\paragraph{Motion Equations in Body-Fixed Representation}

The pivot imposes the geometric constraints%
\begin{equation}
\mathbf{g}\left( C\right) =\mathbf{r}_{0}-\mathbf{R}^{T}\mathbf{r}^{s}=%
\mathbf{0.}  \label{geomContsTop}
\end{equation}%
Time differentiation yields the velocity and acceleration constraints%
\begin{equation}
\left( 
\begin{array}{cc}
\widehat{\mathbf{R}^{T}\mathbf{r}^{s}} & \mathbf{I}%
\end{array}%
\right) \left( 
\begin{array}{c}
\mathbf{\omega } \\ 
\mathbf{v}%
\end{array}%
\right) =\mathbf{JV}=\mathbf{0},\ \ \ \ \ \left( 
\begin{array}{cc}
\widehat{\mathbf{R}^{T}\mathbf{r}^{s}} & \mathbf{I}%
\end{array}%
\right) \left( 
\begin{array}{c}
\dot{\mathbf{\omega }} \\ 
\dot{\mathbf{v}}%
\end{array}%
\right) =-\widehat{\mathbf{\omega }}\widehat{\mathbf{\omega }}\mathbf{r}_{0}+%
\widehat{\mathbf{\omega }}\mathbf{v}  \label{accConstrTop}
\end{equation}%
where the body-fixed twist is denoted $\mathbf{V}=\left( \mathbf{\omega },%
\mathbf{v}\right) $. These together with the body-fixed Newton-Euler
equations w.r.t. to the COM give rise to the system (\ref{dyn4})

\begin{equation}
\left( 
\begin{array}{ccc}
{\Theta }_{0} & \mathbf{0} & -\widehat{\mathbf{r}}_{0} \\ 
\mathbf{0} & m\mathbf{I} & \mathbf{I} \\ 
\widehat{\mathbf{r}}_{0} & \mathbf{I} & \mathbf{0}%
\end{array}%
\right) \left( 
\begin{array}{c}
\dot{\mathbf{\omega }} \\ 
\dot{\mathbf{v}} \\ 
\mathbf{\lambda }%
\end{array}%
\right) =\left( 
\begin{array}{c}
\mathbf{M}-\widehat{\mathbf{\omega }}{\Theta }_{0}\mathbf{\omega } \\ 
\mathbf{F}-m\widehat{\mathbf{\omega }}\mathbf{v} \\ 
-\widehat{\mathbf{\omega }}\widehat{\mathbf{\omega }}\mathbf{r}_{0}+\widehat{%
\mathbf{\omega }}\mathbf{v}%
\end{array}%
\right)  \label{TopEOMbody}
\end{equation}%
where (\ref{geomContsTop}) is assumed satisfied. $\mathbf{F}$ and $\mathbf{M}
$ is the external force and torque, respectively, acting on the COM
represented in the body-fixed frame.

\paragraph{Motion Equations in Hybrid Velocity Representation}

In hybrid velocity representation $\left( \mathbf{\omega },\mathbf{v}%
^{s}\right) $ the geometric constraints (\ref{geomContsTop}) gives rise to
the following velocity and acceleration constraints, respectively,%
\begin{equation}
\left( 
\begin{array}{cc}
\mathbf{R}\widehat{\mathbf{r}}_{0} & \mathbf{I}%
\end{array}%
\right) \left( 
\begin{array}{c}
\mathbf{\omega } \\ 
\mathbf{v}^{s}%
\end{array}%
\right) =\mathbf{0},\ \ \ \ \left( 
\begin{array}{cc}
\mathbf{R}\widehat{\mathbf{r}}_{0} & \mathbf{I}%
\end{array}%
\right) \left( 
\begin{array}{c}
\dot{\mathbf{\omega }} \\ 
\dot{\mathbf{v}}^{s}%
\end{array}%
\right) =\mathbf{R}\widehat{\mathbf{\omega }}\widehat{\mathbf{\omega }}%
\mathbf{r}_{0}.
\end{equation}%
The Newton-Euler equations w.r.t. to COM in hybrid representation yields%
\begin{equation}
\left( 
\begin{array}{ccc}
{\Theta }_{0} & \mathbf{0} & (\mathbf{R}\widehat{r}_{0})^{T} \\ 
\mathbf{0} & m\mathbf{I} & \mathbf{I} \\ 
\mathbf{R}\widehat{\mathbf{r}}_{0} & \mathbf{I} & \mathbf{0}%
\end{array}%
\right) \left( 
\begin{array}{c}
\dot{\mathbf{\omega }} \\ 
\dot{\mathbf{v}}^{s} \\ 
\mathbf{\lambda }%
\end{array}%
\right) =\left( 
\begin{array}{c}
\mathbf{M}-\widehat{\mathbf{\omega }}{\Theta }_{0}\mathbf{\omega } \\ 
\mathbf{F}^{s} \\ 
\mathbf{R}\widehat{\mathbf{\omega }}\widehat{\mathbf{\omega }}\mathbf{r}_{0}%
\end{array}%
\right) .  \label{TopEOMspatial}
\end{equation}%
$\mathbf{F}^{s}$ is the external force acting on the COM represented in the
spatial frame.

\paragraph{Numerical Results}

No gravity or external forces and torques are present, i.e. $\mathbf{M}=%
\mathbf{F}=\mathbf{0}$. The integration step size is $\Delta t=10^{-3}$~s.
In the reference configuration the top has the same orientation as the
inertial frame shown in figure \ref{figHaevyTop}a). The initial angular
velocity was set to $\mathbf{\omega }_{0}=(0,20\,\pi ,10\,\pi )$\thinspace
rad/s so that the top will perform a spatial rotation.

\begin{figure}[h]
a)%
\centerline{
\includegraphics[width=6.7cm]{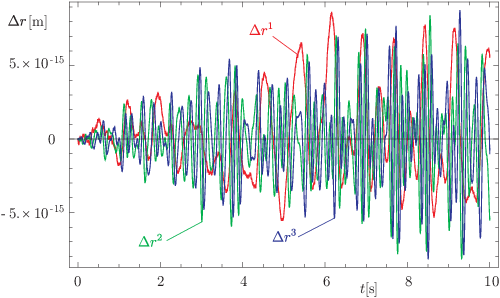}
} b)%
\centerline{
\includegraphics[width=6.7cm]{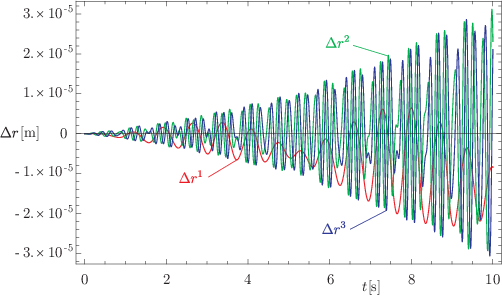}
} \vspace{-4ex}
\caption{Drift of rotation center for integration on a) $SE(3)$, and b)$%
SO(3)\times \mathbb{R}^{3}$. }
\label{figPosErrorSkew}\vspace{-6ex}
\end{figure}

Figure \ref{figPosErrorSkew} shows the deviation of the location of the COM
reference frame for numerical solutions obtained with the $SE\left( 3\right) 
$ and $SO\left( 3\right) \times \mathbb{R}^{3}$ formulation, respectively.
The $SE\left( 3\right) $ integration yields the correct result within the
computation accuracy while the $SO\left( 3\right) \times \mathbb{R}^{3}$
formulation shows significant drifts. \vspace{-3ex}

\begin{figure}[h]
a)%
\centerline{
\includegraphics[width=6.7cm]{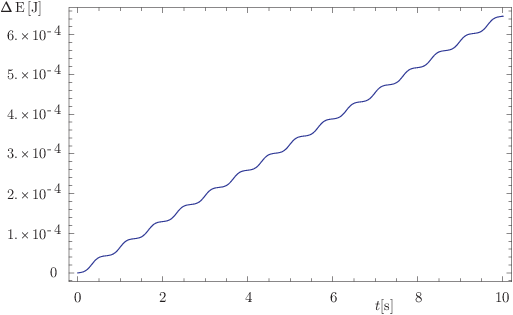}}
b)%
\centerline{\includegraphics[width=6.7cm]{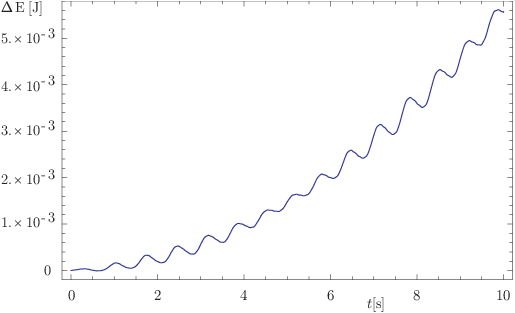}
} \vspace{-2ex} 
\caption{Drift of kinetic energy for a) $SE(3)$ and b) $SO(3)\times \mathbb{R%
}^{3}$ formulation.}
\label{figSE3_EKinError}
\vspace{-3ex}
\end{figure}

A reference trajectory was numerically determined by integrating the dynamic
Euler-equations in quaternion parameterization using a RK4 integration
scheme with variable step size and relative and absolute accuracy goal of $%
10^{-6}$ and $10^{-9}$, respectively. Figure \ref{figSE3_EKinError} shows
the kinetic energy drift for integration on $SE(3)$ relative to the initial
value $T_{0}=13.97$~kJ. This is two orders of magnitude smaller than the
drift of the $SO(3)\times \mathbb{R}^{3}$ formulation (fig. \ref%
{figSE3_EKinError}b)).

\subsection{Spherical Double Pendulum in Gravity Field}

\begin{figure}[h]
\centerline{
\includegraphics[width=5.0cm]{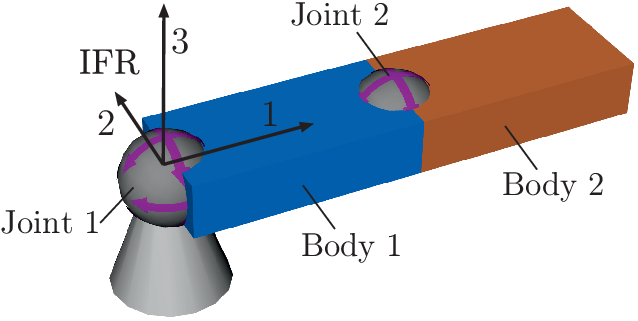}
} 
\caption{Spherical double pendulum.}
\label{figSphericalPendulum}\vspace{-4ex}
\end{figure}

The double pendulum is considered consisting of two slender rigid bodies as
shown in figure \ref{figSphericalPendulum}b). The two bodies are
interconnected and the pendulum as a whole is connected to the ground by
spherical joints. The two links are flat boxes with side length $a,b,c$
along the axes of the COM reference frame. Both have the same dimension with
lengths $a=0.2\,$m$,b=0.1\,$m$,c=0.05\,$m. The configuration of the system
is represented by $C_{1}=\left( \mathbf{R}_{1},\mathbf{r}_{1}^{s}\right) $
and $C_{2}=\left( \mathbf{R}_{2},\mathbf{r}_{2}^{s}\right) $. The two links
are subject to the geometric constraints \vspace{-5ex}

\begin{eqnarray}
\mathbf{g}_{1}\left( C_{1}\right) &=&\mathbf{R}_{1}\mathbf{r}_{0}+\mathbf{r}%
_{1}^{s}=\mathbf{0}  \nonumber \\
\mathbf{g}_{2}\left( C_{1},C_{2}\right) &=&\mathbf{R}_{1}\mathbf{r}_{10}+%
\mathbf{r}_{1}^{s}-\mathbf{R}_{2}\mathbf{r}_{20}-\mathbf{r}_{2}^{s}=\mathbf{0%
}  \label{GeomConstrPend}
\end{eqnarray}
\vspace{-4ex}

where $\mathbf{r}_{i0},i=1,2$ is the position vector from the COM frame on
body $i$ to the spherical joint connecting the two links, expressed in the
COM frame. $\mathbf{r}_{0}$ is the position vector from the COM frame on
body 1 to the spherical joint connecting it to the ground expressed in this
COM frame. Denote with ${\Theta }_{i0}$ the inertia tensor of body $i$
w.r.t. the COM frame, and with $m_{i}$ its mass.

\paragraph{Motion Equations in Body-Fixed Representation}

The velocity and acceleration constraints corresponding to (\ref%
{GeomConstrPend}), in terms of body-fixed twists $\mathbf{V}_{1},\mathbf{V}%
_{2}\in se\left( 3\right) $, are \vspace{-5ex}

\begin{eqnarray}
\mathbf{0} &=&\left( 
\begin{array}{cccc}
\widehat{\mathbf{r}}_{0} & -\mathbf{I} & \mathbf{0} & \mathbf{0} \\ 
\mathbf{R}_{1}\widehat{\mathbf{r}}_{10} & -\mathbf{R}_{1} & -\mathbf{R}_{2}%
\widehat{\mathbf{r}}_{20} & \mathbf{R}_{2}%
\end{array}%
\right) \left( 
\begin{array}{c}
\mathbf{\omega }_{1} \\ 
\mathbf{v}_{1} \\ 
\mathbf{\omega }_{2} \\ 
\mathbf{v}_{2}%
\end{array}%
\right) =\mathbf{JV}  \label{velConstrPendulum} \\
\mathbf{J}\dot{\mathbf{V}} &=&\left( 
\begin{array}{c}
\widehat{\mathbf{\omega }}_{1}\mathbf{v}_{1}\mathbf{+}\widehat{\mathbf{%
\omega }}_{1}\widehat{\mathbf{\omega }}_{1}\mathbf{r}_{0} \\ 
\mathbf{R}_{1}\widehat{\mathbf{\omega }}_{1}\widehat{\mathbf{\omega }}_{1}%
\mathbf{r}_{10}-\mathbf{R}_{2}\widehat{\mathbf{\omega }}_{2}\widehat{\mathbf{%
\omega }}_{2}\mathbf{r}_{20}+\mathbf{R}_{1}\widehat{\mathbf{\omega }}_{1}%
\mathbf{v}_{1}-\mathbf{R}_{2}\widehat{\mathbf{\omega }}_{2}\mathbf{v}_{2}%
\end{array}%
\right) .  \label{accConstrPendulum}
\end{eqnarray}%
\vspace{-3ex}

Using a reference frame at the COM yields \vspace{-4ex}

\begin{equation}
\left( 
\begin{array}{cccccc}
{\Theta }_{10} & \mathbf{0} & \mathbf{0} & \mathbf{0} & -\widehat{\mathbf{r}}%
_{0} & -\widehat{\mathbf{r}}_{10}\mathbf{R}_{1}^{T} \\ 
\mathbf{0} & m_{1}\mathbf{I} & \mathbf{0} & \mathbf{0} & -\mathbf{I} & -%
\mathbf{R}_{1}^{T} \\ 
\mathbf{0} & \mathbf{0} & {\Theta }_{20} & \mathbf{0} & \mathbf{0} & 
\widehat{\mathbf{r}}_{20}\mathbf{R}_{2}^{T} \\ 
\mathbf{0} & \mathbf{0} & \mathbf{0} & m_{2}\mathbf{I} & \mathbf{0} & 
\mathbf{R}_{2}^{T} \\ 
\widehat{\mathbf{r}}_{0} & -\mathbf{I} & \mathbf{0} & \mathbf{0} & \mathbf{0}
& \mathbf{0} \\ 
\mathbf{R}_{1}\widehat{\mathbf{r}}_{10} & -\mathbf{R}_{1} & -\mathbf{R}_{2}%
\widehat{\mathbf{r}}_{20} & \mathbf{R}_{2} & \mathbf{0} & \mathbf{0}%
\end{array}%
\right) 
\hspace{-1.2ex}%
\left( 
\begin{array}{c}
\dot{\mathbf{\omega }}_{1} \\ 
\dot{\mathbf{v}}_{1} \\ 
\dot{\mathbf{\omega }}_{2} \\ 
\dot{\mathbf{v}}_{2} \\ 
\mathbf{\lambda }_{1} \\ 
\mathbf{\lambda }_{2}%
\end{array}%
\right) 
\hspace{-1.2ex}%
=%
\hspace{-1ex}%
\left( 
\hspace{-1ex}%
\begin{array}{c}
-\widehat{\mathbf{\omega }}_{1}{\Theta }_{10}\mathbf{\omega }_{1} \\ 
\mathbf{F}_{1}-m_{1}\widehat{\mathbf{\omega }}_{1}\mathbf{v}_{1} \\ 
-\widehat{\mathbf{\omega }}_{2}{\Theta }_{20}\mathbf{\omega }_{2} \\ 
\mathbf{F}_{2}-m_{2}\widehat{\mathbf{\omega }}_{2}\mathbf{v}_{2} \\ 
\ast \\ 
\ast \ast%
\end{array}%
\right)
\end{equation}%
where $\ast $ and $\ast \ast $ are the terms in the two rows of the right
hand side of (\ref{accConstrPendulum}). Since only gravity forces act on the
systems the body-fixed forces are $\mathbf{F}_{i}=\mathbf{R}_{i}^{T}\mathbf{g%
}^{s}$, where $\mathbf{g}^{s}=\left( 0,0,-g\right) $ is the gravity vector
w.r.t. space-fixed frame. The Lagrange multiplier $\mathbf{\lambda }_{i}\in 
\mathbb{R}^{3}$ is the reaction force in joint $i$.

\paragraph{Motion Equations in Hybrid Velocity Representation}

With hybrid velocities $\mathbf{V}_{i}=\left( \mathbf{\omega }_{i},\mathbf{v}%
_{i}^{s}\right) $ the velocity and acceleration constraints are%
\vspace{-6ex}%

\begin{eqnarray}
\mathbf{0} &=&\left( 
\begin{array}{cccc}
\mathbf{R}_{1}\widehat{\mathbf{r}}_{0} & -\mathbf{I} & \mathbf{0} & \mathbf{0%
} \\ 
\mathbf{R}_{1}\widehat{\mathbf{r}}_{10} & -\mathbf{I} & -\mathbf{R}_{2}%
\widehat{\mathbf{r}}_{20} & \mathbf{I}%
\end{array}%
\right) \left( 
\begin{array}{c}
\mathbf{\omega }_{1} \\ 
\mathbf{v}_{1} \\ 
\mathbf{\omega }_{2} \\ 
\mathbf{v}_{2}%
\end{array}%
\right) =\mathbf{JV}  \label{velConstrPendulumHyb} \\
\mathbf{J}\dot{\mathbf{V}} &=&\left( 
\begin{array}{c}
\mathbf{R}_{1}\widehat{\mathbf{\omega }}_{1}\widehat{\mathbf{\omega }}_{1}%
\mathbf{r}_{0} \\ 
\mathbf{R}_{1}\widehat{\mathbf{\omega }}_{1}\widehat{\mathbf{\omega }}_{1}%
\mathbf{r}_{10}-\mathbf{R}_{2}\widehat{\mathbf{\omega }}_{2}\widehat{\mathbf{%
\omega }}_{2}\mathbf{r}_{20}%
\end{array}%
\right) .  \label{accConstrPendulumHyb}
\end{eqnarray}%
\vspace{-4ex}

The index 1 motion equations are, with force vectors $\mathbf{F}_{i}^{s}=%
\mathbf{g}^{s}$, \vspace{-3ex}

\[
\left( 
\begin{array}{cccccc}
{\Theta }_{10} & \mathbf{0} & \mathbf{0} & \mathbf{0} & -\widehat{\mathbf{r}}%
_{0}\mathbf{R}_{1}^{T} & -\widehat{\mathbf{r}}_{10}\mathbf{R}_{1}^{T} \\ 
\mathbf{0} & m_{1}\mathbf{I} & \mathbf{0} & \mathbf{0} & -\mathbf{I} & -%
\mathbf{I} \\ 
\mathbf{0} & \mathbf{0} & {\Theta }_{20} & \mathbf{0} & \mathbf{0} & 
\widehat{\mathbf{r}}_{20}\mathbf{R}_{2}^{T} \\ 
\mathbf{0} & \mathbf{0} & \mathbf{0} & m_{2}\mathbf{I} & \mathbf{0} & 
\mathbf{I} \\ 
\mathbf{R}_{1}\widehat{\mathbf{r}}_{0} & -\mathbf{I} & \mathbf{0} & \mathbf{0%
} & \mathbf{0} & \mathbf{0} \\ 
\mathbf{R}_{1}\widehat{\mathbf{r}}_{10} & -\mathbf{I} & -\mathbf{R}_{2}%
\widehat{\mathbf{r}}_{20} & \mathbf{I} & \mathbf{0} & \mathbf{0}%
\end{array}%
\right) 
\hspace{-1.2ex}%
\left( 
\begin{array}{c}
\dot{\mathbf{\omega }}_{1} \\ 
\dot{\mathbf{v}}_{1}^{s} \\ 
\dot{\mathbf{\omega }}_{2} \\ 
\dot{\mathbf{v}}_{2}^{s} \\ 
\mathbf{\lambda }_{1} \\ 
\mathbf{\lambda }_{2}%
\end{array}%
\right) 
\hspace{-1.2ex}%
=%
\hspace{-1ex}%
\left( 
\hspace{-1ex}%
\begin{array}{c}
-\widehat{\mathbf{\omega }}_{1}{\Theta }_{10}\mathbf{\omega }_{1} \\ 
\mathbf{F}_{1}^{s} \\ 
-\widehat{\mathbf{\omega }}_{2}{\Theta }_{20}\mathbf{\omega }_{2} \\ 
\mathbf{F}_{2}^{s} \\ 
\ast \\ 
\ast \ast%
\end{array}%
\right) 
\]%
\vspace{-3ex}

The terms $\ast $ and $\ast \ast $ in the two rows of the right hand side of
(\ref{accConstrPendulumHyb}).

\paragraph{Numerical Results}

In the initial configuration $C_{1}\left( 0\right) =\left( \mathbf{I},\left(
a_{1}/2,0,0\right) \right) $, $C_{2}\left( 0\right) =\left( \mathbf{I}%
,\left( a_{1}+a_{2}/2,0,0\right) \right) $ the pendulum is aligned along the
space-fixed x-axis as shown in figure \ref{figSphericalPendulum}. The
initial velocities are set to $\mathbf{\omega }_{10}=(10,0,0)$\thinspace
rad/s and $\mathbf{\omega }_{20}=(10\,\pi ,10\,\pi ,20\,\pi )$\thinspace
rad/s. The pendulum is moving in the gravity field. A time step size of $%
\Delta t=10^{-3}$s is used in the MK method (\ref{MK1}). Figure \ref%
{DoublePendPosError1} and \ref{DoublePendPosError2} reveal a generally
better performance of the $SE\left( 3\right) $ formulation. It also reveals
that the constraint violation is more pronounced when it is due to the
relative motion of two screw motions, which explains that the violation in
joint 2 is higher than that in joint 1. The total energy drift is comparable
for both formulations (figure \ref{DoublePendEnergyDrift}). \vspace{-4ex}

\begin{figure}[th]
a) 
\centerline{
\includegraphics[width=6.7cm]{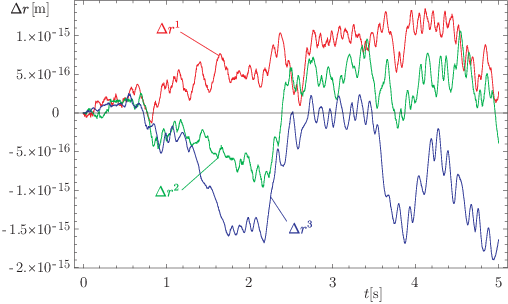}}
b) 
\centerline{
\includegraphics[width=6.7cm]{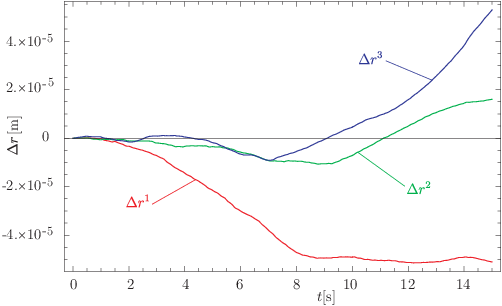}
} \vspace{-2ex}
\caption{Violation of geometric constraints of joint 1 when integrating the $%
SE(3)$ (a), and $SO(3)\times \mathbb{R}^{3}$ (b) formulation. \protect%
\vspace{-2ex}}
\label{DoublePendPosError1}
\end{figure}
\vspace{-4ex}

\begin{figure}[b]
a)\centerline{%
\includegraphics[width=6.7cm]{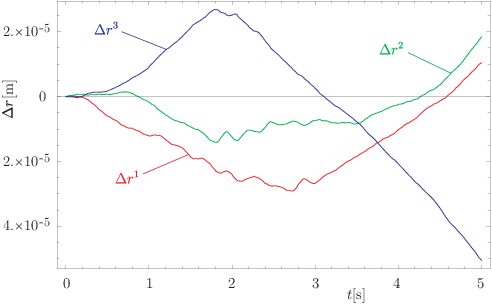}}
b)\centerline{%
\includegraphics[width=6.7cm]{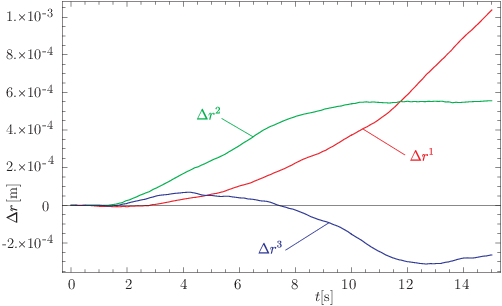}}
\caption{Violation of geometric constraints of joint 2 when integrating a)
the $SE(3)$, and b) $SO(3)\times \mathbb{R}^{3}$ formulation. \protect%
\vspace{-2ex}}
\label{DoublePendPosError2}
\end{figure}
\vspace{-5ex}

\begin{figure}[h]
a)\centerline{%
\includegraphics[width=6.7cm]{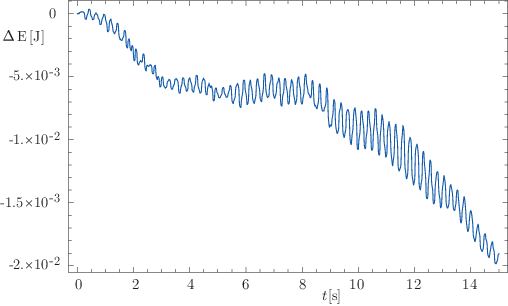}}
b)\centerline{%
\includegraphics[width=6.7cm]{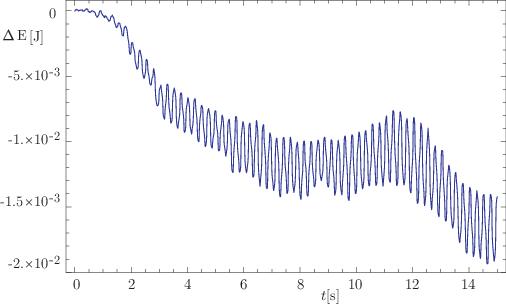}}
\vspace{-2ex}
\caption{Drift of total energy for a) the $SE(3)$, and b) $SO(3)\times 
\mathbb{R}^{3}$. \protect\vspace{-0.5ex}}
\label{DoublePendEnergyDrift}
\end{figure}

\subsection{Interconnected Floating Bodies}

\begin{figure}[h]
\centerline{
\includegraphics[width=5cm]{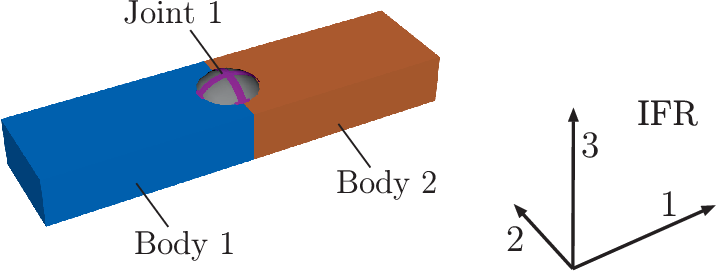}
} \vspace{1.5ex} \vspace{-3ex}
\caption{Two floating bodies connected by spherical joint.}
\label{figFloating}
\end{figure}

Consider the situation in figure \ref{figFloating} where the two bodies in
the above spherical double pendulum mutually connected but are not connected
to the ground. That is, the two bodies, connected by one spherical joint,
are free floating, and no gravity is assumed. The motion equations are the
same as those of the spherical double pendulum above except that the first
constraint in (\ref{velConstrPendulum}),(\ref{accConstrPendulum}), and (\ref%
{velConstrPendulumHyb}),(\ref{accConstrPendulumHyb}), respectively, are
removed. 

Again the initial configuration is such that the bodies are horizontally
alligned. The initial angular velocities are $\mathbf{\omega }%
_{10}=(0,0,-10) $\thinspace rad/s and $\mathbf{\omega }_{20}=(1,-1,2\,\pi )$%
\thinspace rad/s.

\begin{figure}[h]
a)%
\centerline{
\includegraphics[width=6.7cm]{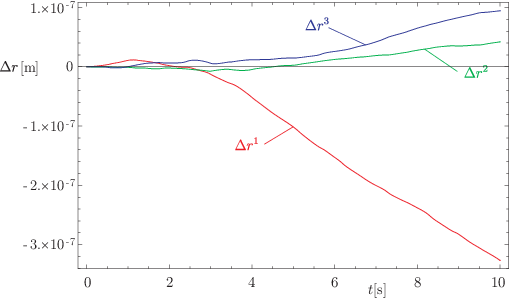}
} b)%
\centerline{
\includegraphics[width=6.7cm]{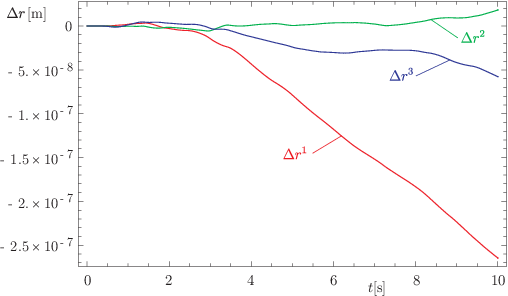}
}  \vspace{-5ex}
\caption{Drift of rotation center for integration on a) $SE(3)$, and b) $%
SO(3)\times \mathbb{R}^{3}$. }
\label{figPosErrorFloatingDouble}
\end{figure}

The numerical results obtained by the Lie group integration schemes are
shown in figure \ref{figPosErrorFloatingDouble}. Apparently for this example
both formulations perform similarly regarding constraint satisfaction. Also
the kinetic energy drifts are similar as shown in \ref%
{figEKinErrorFloatingDouble}.

This example differs from the previous ones in that the constraint violation
is caused by the spatial motion of two interconnected bodies rather then by
the motion of one body that is connected to the ground.
\vspace{-5ex}

\begin{figure}[h]\vspace{-1ex}
a)%
\centerline{
\includegraphics[width=6.7cm]{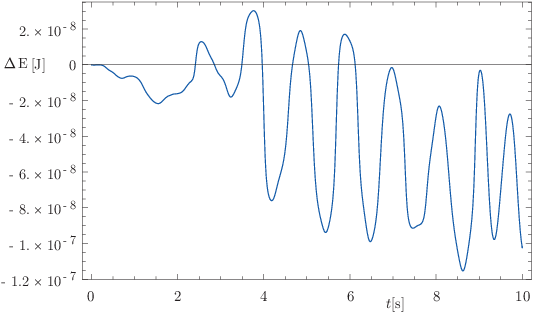}}
b)%
\centerline{\includegraphics[width=6.7cm]{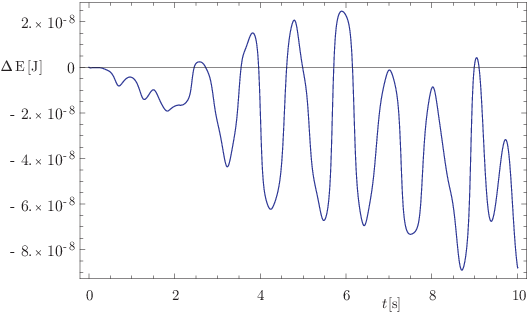}
} \vspace{-2ex} \vspace{-3ex}
\caption{Drift of kinetic energy for a) $SE(3)$ and b) $SO(3)\times \mathbb{R%
}^{3}$.}
\label{figEKinErrorFloatingDouble}
\end{figure}

\vspace{-5ex}
\subsection{Closed Loop Spherical 3-Bar Linkage}

\begin{figure}[h]
\centerline{
\includegraphics[width=4.8cm]{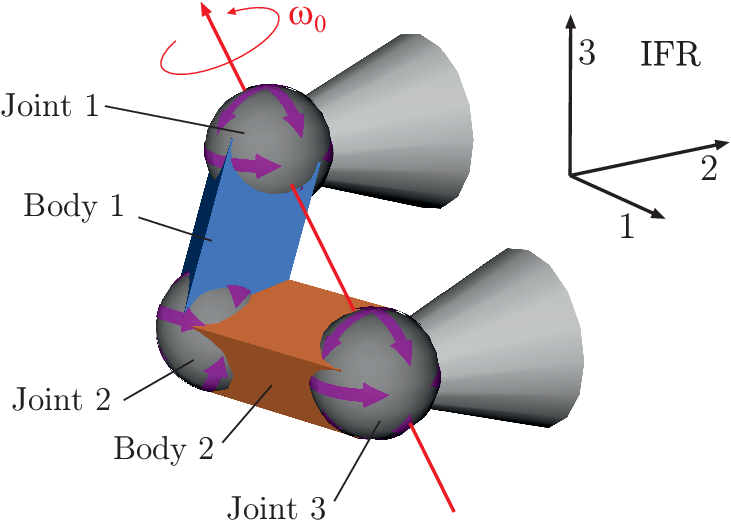}
} \vspace{-1ex}
\caption{Closed loop 3-bar linkage.}
\label{fig3Bar}
\end{figure}

\begin{figure}[t]
a) 
\centerline{
\includegraphics[width=6.4cm]{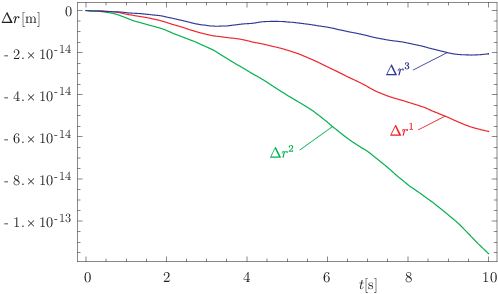}}
b) 
\centerline{
\includegraphics[width=6.4cm]{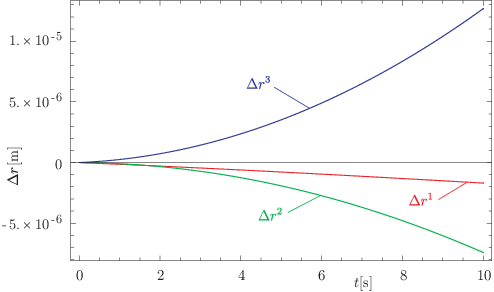}
} \vspace{-4ex}
\caption{Violation of geometric constraints of joint 1 when integrating the $%
SE(3)$ (a), and $SO(3)\times \mathbb{R}^{3}$ (b) formulation. \protect%
\vspace{-5ex}}
\label{ClosedLoopPosError1}
\end{figure}

As final example the closed-loop 3-bar mechanism in figure \ref{fig3Bar}
is considered. Due to space limitations the motion equation are not given. The results in figure \ref{ClosedLoopPosError1}-\ref%
{ClosedLoopPosError3} confirm that again that the constraints are perfectly
satisfied if only the motion of one body is to be estimated, i.e. the base
joints. The constraint violation of joint 2 connecting the two bodies is
similar for both formulations. Also the energy drift of both formulations is
similar (fig. \ref{ClosedLoopEnergyDrift}). \vspace{-3ex}

\begin{figure}[b]
a) 
\centerline{
\includegraphics[width=6.4cm]{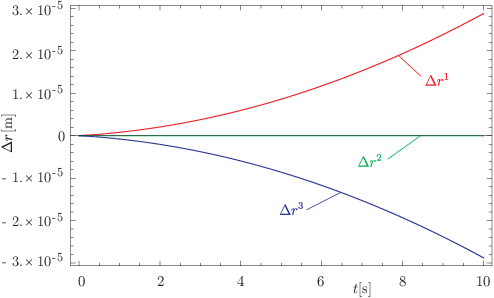}}
b) 
\centerline{
\includegraphics[width=6.4cm]{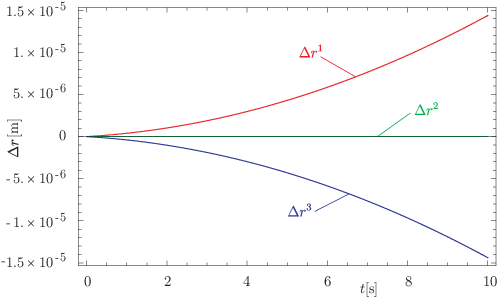}
} \vspace{-4ex}
\caption{Violation of geometric constraints of joint 2 when integrating the $%
SE(3)$ (a), and $SO(3)\times \mathbb{R}^{3}$ (b) formulation. \protect%
\vspace{-4ex}}
\label{ClosedLoopPosError2}
\end{figure}
\begin{figure}[h]
a) 
\centerline{
\includegraphics[width=6.7cm]{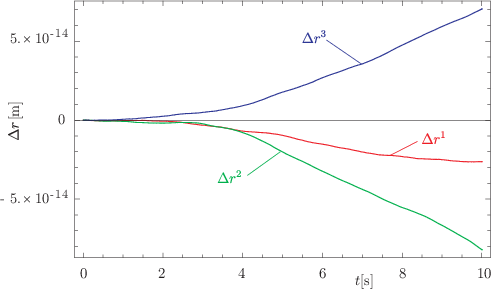}}
b) 
\centerline{
\includegraphics[width=6.7cm]{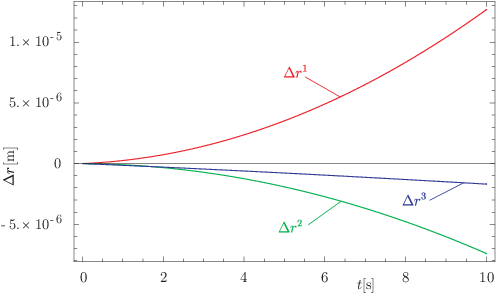}
} \vspace{-2ex}
\caption{Violation of geometric constraints of joint 3 when integrating the $%
SE(3)$ (a), and $SO(3)\times \mathbb{R}^{3}$ (b) formulation. \protect%
\vspace{-2ex}}
\label{ClosedLoopPosError3}
\end{figure}

\begin{figure}[h]
a)\centerline{%
\includegraphics[width=6.7cm]{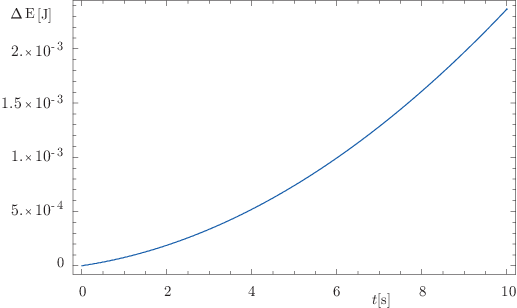}}
b)\centerline{%
\includegraphics[width=6.7cm]{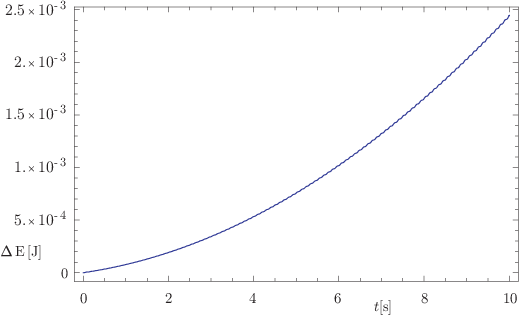}}
\vspace{-2ex}
\caption{Drift of total energy for a) the $SE(3)$, and b) $SO(3)\times 
\mathbb{R}^{3}$. \protect\vspace{-2ex}}
\label{ClosedLoopEnergyDrift}
\end{figure}

\section{Discussion and Conclusion}

A generic rigid body motion is a screw motion, and the corresponding
velocity is a proper twist. Consequently the finite motion of a body is best
approximated from its instantaneous twist as finite screw motion. As the
actual cause of the body twist is irrelevant this applies equally to a free
floating body as well as to a holonomically and scleronomically constrained
body, in particular when jointly connected to the ground. From a
computational point of view the numerical representation of screw motions is
crucial. Generally, the Lie group $SE\left( 3\right) $ represents rigid body
motions and is the proper rigid body configuration space allowing for the
reconstruction of finite motions from velocities. Nevertheless, frequently
the direct product group $SO\left( 3\right) \times \mathbb{R}^{3}$ is
employed as configuration space. The latter does allow for representing
rigid body configurations but not motions, however. The question arises
whether this is significant for the numerical simulation. Lie group
integration schemes take explicitly into account the configuration space
manifold, so that the numerical update within the integration scheme is
different for the two configuration spaces. With the above said it is clear
that the actual performance of the integration scheme depends on the
underlying configuration space Lie group, and it may be conjectured that the 
$SE\left( 3\right) $ formulation outperforms the $SO\left( 3\right) \times 
\mathbb{R}^{3}$ formulation.

In this paper the accuracy of numerical Munthe-Kaas integration schemes
applied to the two formulations is investigated. The results confirm that
the $SE\left( 3\right) $ formulation yields the best overall performance. In
particular, since rigid body motions are described relative to an inertial
reference frame, the $SE\left( 3\right) $ formulation allows for perfect
reconstruction of the motion of a constrained rigid body if its motion is
constrained to a proper motion subgroup. This is the case when a rigid body
is connected to the ground (or to another slowly moving body) by a lower
pair. In the general case of rigid bodies connected by lower pairs their
relative motion belongs to a motion subgroup but their absolute motion does
not. In these cases also the reconstruction as screw motion of the
individual bodies from first-order motions cannot capture the
interdependence of their finite motions, and both formulations perform
equally.

In summary the $SE\left( 3\right) $ update scheme performs generally as good
as the $SO\left( 3\right) \times \mathbb{R}^{3}$ formulation while it
outperforms the latter when bodies are constrained to a stationary body, in
particular the ground.

This advantage is owed to an increase in cmplexitiy since the $SE\left(
3\right) $ update scheme is computationally more complex than the standard
direct product $SO\left( 3\right) \times \mathbb{R}^{3}$ update. In order to
minimize the complexity of the integration scheme, for a given problem, the
configuration space can be introduced so that the $SE\left( 3\right) $ is
used as configuration space where appropriate and $SO\left( 3\right) \times 
\mathbb{R}^{3}$ otherwise. Such tailored designation of configuration spaces
will be part of furure work aiming on integration schemes that minimize
constraint violation in numerical MBS models. This paper addressed the
performance of MK integration methods. To complete this study the effect of
different configuration space for other Lie group integration schemes will
be investigated in future.


\begin{thebibliography}{99}
\bibitem{Blajer2011} W. Blajer: Methods for constraint violation suppression
in the numerical simulation of constrained multibody systems -- A
comparative study, Comput. Methods Appl. Mech. Engrg. Vol. 200, 2011, pp.
1568--1576

\bibitem{BottemaRoth1979} A. Bottema, B. Roth: Theoretical Kinematics,
North-Holland, 1979

\bibitem{BruelsCardona2010} O. Br\"{u}ls, A. Cardona: On the Use of Lie
Group Time Integrators in Multibody Dynamics, J. Comput. Nonlinear Dynam.
5(3) 2010

\bibitem{BruelsCardonaArnold2012} O. Br\"{u}ls, A. Cardona, M. Arnold: Lie
group generalized-alpha time integration of constrained flexible multibody
systems, Mech. Mach. Theory 48, 2012, 121-137

\bibitem{Bruyninckx1996} H. Bruyninckx, J. De Schutter: Symbolic
differentiation of the velocity mapping for a serial kinematic chain, Mech.
Mach. Theory 31(2) 1996, 135-148

\bibitem{BulloMurray1995} F. Bullo, R.M. Murray: \textit{Proportional
Derivative (PD) Control on the Euclidean group}, CDS technical report
95-010, 1995

\bibitem{CelledoniOwren1999} E. Celledoni, B. Owren:\textit{\ Lie Group
Methods for Rigid Body Dynamics and Time Integration on Manifolds}, Computer
Methods in Applied Mechanics and Engineering 19, 1999, 421-438

\bibitem{Hairer2006} E. Hairer, C. Lubich and G. Wanner: Geometric Numerical
Integration, Springer, 2006

\bibitem{Krysl} P. Krysl, L. Endres: Explicit Newmark/Verlet algorithm for
time integration of the rotational dynamics of rigid bodies, Int. J. Numer.
Meth. Engng. 62, 2005, 2154-2177

\bibitem{muntekaas1998} H. Munthe-Kaas, \textit{Runge Kutta methods on Lie
groups}, BIT, 38(1) 1998, 92-111

\bibitem{muntekaas1999} H. Munthe-Kaas, \textit{High order Runge-Kutta
methods on manifolds}, Appl. Numer. Math. 29, 1999, 115-127

\bibitem{murray} R.M Murray, Z. Li, S. S. Sastry: \textit{A mathematical
Introduction to robotic Manipulation}, CRC Press, 1993

\bibitem{SeligBook} J.M. Selig, \textit{Geometrical Methods in Robotics},
Springer, New York, 1996.

\bibitem{Terze2011} Z. Terze, A. M\"{u}ller, D. Zlatar: Lie-Group
integration method for constrained multibody systems in stabilized index 1
form, Multibody System Dynamics, submitted, 2011
\end{thebibliography}
\end{document}